\documentclass[11pt]{amsart}

\usepackage{amsmath}

\usepackage{amssymb}

\usepackage{graphicx}

\newtheorem{thm}{Theorem}

\newtheorem{lemma}[thm]{Lemma}

\theoremstyle{definition}

\theoremstyle{remark}
\newtheorem{remark}[thm]{Remark}

\def\C{\mathbb{C}}
\def\R{\mathbb{R}}

\def\H{\mathbb{H}}

\def\sl{\mathfrak{sl}}

\def\SO{\mathrm{SO}}
\def\SL{\mathrm{SL}}
\def\GL{\mathrm{GL}}

\def\det{\mathrm{det\,}}
\def\dim{\mathrm{dim\,}}

\begin{document}

\title[Closed trajectories on null curves]{Closed trajectories
 of a particle model \\
on null curves in anti-de Sitter 3-space}

\author{Emilio Musso}
\address{(E. Musso) Dipartimento di Matematica Pura ed Applicata,
Universit\`a degli Studi dell'Aquila, Via Vetoio, I-67010
Coppito (L'Aquila), Italy} \email{musso@univaq.it}

\author{Lorenzo Nicolodi}
\address{(L. Nicolodi) Di\-par\-ti\-men\-to di Ma\-te\-ma\-ti\-ca,
Uni\-ver\-si\-t\`a degli Studi di Parma, Viale G. P. Usberti 53/A,
I-43100 Parma, Italy} \email{lorenzo.nicolodi@unipr.it}

\thanks{Authors partially supported by MIUR projects:
\textit{Metriche riemanniane e variet\`a differenziali} (E.M.);
\textit{Propriet\`a geometriche delle variet\`a reali e complesse} (L.N.);
and by the GNSAGA of INDAM}

\subjclass[2000]{58E10; 49F05}

%\date{Version of September 12, 2007}

%\dedicatory{}

\keywords{Null curves, closed trajectories, anti-de Sitter 3-space}

\begin{abstract}
We study the existence of closed trajectories of a particle model on
null curves in anti-de Sitter 3-space defined by a functional which
is linear in the curvature of the particle path. Explicit
expressions for the trajectories are found and the existence of
infinitely many closed trajectories is proved.
\end{abstract}

\maketitle

\section{Introduction}\label{s:intro}

In this paper we study null curves in anti-de Sitter 3-space
which are critical points for the functional
\begin{equation}\label{functional}
 \mathcal{L}(\gamma) = \int_\gamma{(m + k_\gamma)}ds, \quad m\in \R,
  \end{equation}
where $s$ is the pseudo-arc parameter which normalizes the
derivative of the tangent vector field of $\gamma$ and $k_\gamma(s)$
is a curvature function that, in general, uniquely determines
$\gamma$ up to Lorentz transformations (\cite{FGL-PhysLettB},
\cite{Inoguchi}). The functional \eqref{functional} is invariant
under the group $\SL(2,\R)\times\SL(2,\R)$, which doubly covers the
identity component of the group of Lorentz transformations.
Motivations for this study are provided by optimal control theory
and especially by the recent interest in certain particle models on
null curves in Lorentzian 3-space forms associated with action
integrals of the type above (\cite{NMMK-NucPhysB},
\cite{FGL-PhysLettB}). Yet another motivation is given by surface
geometry; if we take $\SL(2,\R)$ as a model for anti-de Sitter
3-space, then a null curve in anti-de Sitter 3-space, as real form
of a holomorphic null curve in $\SL(2,\C)$, is related to the theory
of constant mean curvature one (cmc-1) surfaces and flat fronts in
hyperbolic 3-space (\cite{Br1987}, \cite{UY-Crelle}, \cite{GMM-MA}).
In perspective, one would like to understand the class of cmc-1
surfaces generated by the critical points of \eqref{functional}.

\vskip0.1cm

The purpose of this article is to investigate the global behavior of
extremal trajectories of the functional \eqref{functional}. We will
find explicit expressions for the extremal trajectories and then
establish the existence of infinitely many closed ones. This result
is related to the presence of maximal compact abelian subgroups in
the isometry group $\SL(2,\R)\times\SL(2,\R)$ of anti-de Sitter
3-space. For a discussion of extremal trajectories in the other
Lorentzian space forms we refer to \cite{GM}, \cite{MN-SICON}.

\vskip0.1cm

The Euler--Lagrange equation associated to \eqref{functional} yields
that the curvature of extremal trajectories is either a constant,
in which case we have null helices, or an elliptic
function\footnote{Possibly a
degenerate one, i.e., an hyperbolic, trigonometric, or rational
function.} of the pseudo-arc parameter. In this case,
extremal trajectories are governed by a second order linear ODE
with doubly periodic coefficients.
%of Lam\'e type.
By classical results of Picard \cite{Pi} in the Fuchsian theory of
linear ODEs, the solution curves are then expressible in terms of
the Weierstrass $\wp$, $\sigma$ and $\zeta$ functions (cf. Theorem
\ref{A}). The explicit integration of extremal trajectories amounts
to the integration of a linearizable flow on $K_1\times K_2$, where
$K_1$ and $K_2$ are 1-dimensional abelian subgroups of $\SL(2,\R)$.
In particular, if $K_1=K_2=\SO(2)$, the integration amounts to
solving a linearizable first-order ODE on a 2-dimensional torus.
This setting strongly suggests the possibility of periodic
solutions. That this is indeed the case is established in Theorem
\ref{B} where the existence of countably many periodic trajectories
is proved by studying the map of periods. The proof relies on
computations made with \textsc{Mathematica}.
%(cf. Figure 2).

\vskip0.1cm

The paper is organized as follows. Section \ref{s:pre} contains some
background material. Section \ref{s:integration} provides the
explicit integration of extremal trajectories. Section
\ref{s:periodic-solutions} discusses periodic trajectories and
proves the existence of infinitely many of them. The Appendix
outlines the derivation of the Euler-Lagrange equation associated
with \eqref{functional} via the Griffiths formalism \cite{Gr}.

\section{Preliminaries}\label{s:pre}

Anti-de Sitter 3-space, $\H^3_1$, can be viewed as the special
linear group $\SL(2,\R)$ endowed with the bi-invariant Lorentz
metric of constant sectional curvature $-1$ defined by the quadratic
form
\begin{equation}
 q(X) = (x^1_1)^2 + x^1_2x^2_1 = - \det X,
    \end{equation}
for each $X=(x^i_j) \in \sl(2,\R)$. The group $G = \SL(2,\R) \times
\SL(2,\R)$ acts transitively by isometries on $\H^3_1\cong
\SL(2,\R)$ via the action
\[
 (A,B)\cdot x = AxB^{-1}.
  \]
The stability subgroup at the identity $I_{2}$ is the diagonal group
\[
   \Delta=\{(A,A) \, | \, A\in \SL(2,\R)\}
     \]
and $\H^3_1$ may be described as a Lorentzian symmetric space
\[
 \H^3_1 \cong G/\Delta.
  \]
The projection
\[
 \pi : G \ni (A,B) \mapsto AB^{-1}\in \SL(2,\R)
  \]
makes $G$ into a principal bundle with structure group $\Delta$.

\vskip0.2cm

Let $I\subset \R$ be any open interval of real numbers. A smooth
parametrized curve $\gamma : I  \to \SL(2,\R)$ is \textit{null}, or
\textit{lightlike}, if $\det(\gamma^{-1}\gamma')$ vanishes
identically. If $\gamma$ has no flex points,\footnote{i.e.,
$\gamma'(t)$ and $\gamma''(t)$ are linearly independent, for each
$t\in I$, where $\gamma''$ is the covariant derivative of $\gamma'$
along the curve.} there exists a canonical lift
\[
 \Gamma = (\Gamma_+, \Gamma_-) : I \to G
  \]
such that
\begin{equation}\label{MC}
 \Gamma_+^{-1}d\Gamma_+ =
  \left(
   \begin{array}{cc}
     0&1\\
      k+1&0\\
       \end{array}\right)\omega,
\quad
\Gamma_-^{-1}d\Gamma_- =
 \left(
  \begin{array}{cc}
    0&1\\
     k-1&0\\
      \end{array}\right)\omega,
       \end{equation}
where $\omega$ is a nowhere vanishing 1-form, the \textit{canonical
arc element}, and $k : I \to \R$ is a smooth function, the
\textit{curvature function}. We call $\Gamma$ the \textit{spinor
frame} field along $\gamma$ and its components $\Gamma_+$ and
$\Gamma_-$ the positive and negative spinor frame, respectively. The
spinor frame $\Gamma$ is essentially unique, in the sense that $\pm
\Gamma$ are the only lifts satisfying \eqref{MC}. Throughout the
paper we will consider null curves without flex points and
parametrized by the natural parameter, i.e., $\omega = ds$ (cf.
\cite{FGL-PhysLettB}, \cite{BFJL-CQG}).

Conversely, for a smooth function $k : I \to \R$, let $H_{\pm}(k) :
I \to \sl(2,\R)$ be
\begin{equation}\label{hamiltoniani}
 H_+(k) = \left(
  \begin{array}{cc}
    0&1\\
     k+1&0\\
      \end{array}\right),\quad
  H_-(k) = \left(
   \begin{array}{cc}
    0&1\\
     k-1&0\\
      \end{array}\right).
       \end{equation}
By solving a linear system of ODEs, there exists a unique (up to
left multiplication)
\[
 \Gamma = (\Gamma_+, \Gamma_-) : I \to G
  \]
such that
\[
 \Gamma_+^{-1}{\Gamma}'_+ = H_+(k), \quad \Gamma_-^{-1}{\Gamma}'_- = H_-(k).
  \]
In particular, $\gamma = \Gamma_+\Gamma_-^{-1} : I \to \SL(2,\R)$ is
a null curve without flex points, parametrized by the natural
parameter and with curvature function $k$.\footnote{The curve
$\gamma$ is uniquely defined up to orientation and time-orientation
preserving isometries.}

In this context, two null curves $\gamma : I \to \SL(2,\R)$,
$\tilde{\gamma} : \tilde{I} \to \SL(2,\R)$ are said
\textit{equivalent} if there exist $c\in\R$ and $A, B \in \SL(2,\R)$
such that $\tilde{\gamma}(s) = A\gamma(s+c)B^{-1}$, for all $s\in
\tilde{I}$.

\vskip0.2cm

\textbf{Definition.} An \textit{extremal trajectory}
(or simply, a \textit{trajectory}) in $\H^3_1$ is a null
curve with non-constat curvature which is a critical point of the
action functional
\begin{equation}
 \mathcal{L}_m(\gamma) = \int_\gamma{(m + k_\gamma)\omega}
   \end{equation}
under compactly supported variations, where the Lagrange multiplier
$m$ is a real constant.

\vskip0.2cm

The Euler-Lagrange equation associated to $\mathcal{L}_m$ is
computed to be
\[
 {k}''' -6k{k}' + 2m{k}' = 0
  \]
(cf. \cite{FGL-PhysLettB}; see also the Appendix for a different way
of deriving this equation). This may be thought of as the intrinsic
equation of a trajectory. If we let $h := \frac{1}{2}(k -
\frac{m}{3})$ be the \textit{reduced curvature}, then $h$ satisfies
\begin{equation}\label{eq-potenziale}
 ({h}')^2 = 4h^3 -g_2h -g_3,
  \end{equation}
for real constants $g_2$ and $g_3$. Hence $h$ is expressed by the
real values of either a Weierstrass $\wp$-function with invariants
$g_2$, $g_3$, or one of its degenerate forms.

\vskip0.2cm

We call a solution to \eqref{eq-potenziale} a \textit{potential}
with analytic invariants $g_2$, $g_3$. Two potentials are considered
equivalent if they differ by a reparametrization of the form $s
\mapsto s+c$, where $c$ is a constant.\footnote{When invariants
$g_2$ and $g_3$ are given, such that $27g_3^2 \neq g_2^3$, the
general solution of the differential equation $(\frac{dy}{dz})^2 =
4y^3 -g_2y -g_3$ can be written in the form $\wp(z
+\alpha;g_2,g_3)$, where $\alpha$ is a constant of integration.} For
real $g_2$ and $g_3$, let $\Delta(g_2,g_3) = 27g_3^2 -g_2^3$ be the
discriminant of the cubic polynomial
\[
 P(t;g_2,g_3) = 4t^3 -g_2t -g_3.
  \]
The study of the real values of the Weierstrass $\wp$-function with
real invariants $g_2$, $g_3$ (and its degenerate forms) leads to
primitive half-periods $\omega_1$, $\omega_3$ such that:
\begin{itemize}

\item $\Delta(g_2,g_3)<0$: $\omega_1>0$, $\omega_3=i\nu\omega_1$, $\nu>0$.

\item $\Delta(g_2,g_3)>0$: $\omega_1>0$, $\omega_3=\frac{1}{2}(1+i\nu)\omega_1$,
$\nu>0$.

\item $\Delta(g_2,g_3)=0$ and $g_3 > 0$: $\omega_1>0$, $\omega_3 =+i\infty$.

\item $\Delta(g_2,g_3)=0$ and $g_3 < 0$: $\omega_1=+\infty$, $-i\omega_3>0$.

\item $g_2=g_3=0$: $\omega_1=+\infty$, $\omega_3=+i\infty$.

\end{itemize}
Accordingly, denoting by $\mathcal{D}(g_2,g_3)$ the fundamental
period-parallelogram spanned by $2\omega_1$ and $2\omega_3$, the
only possible cases for the potential function $h : I \to \R$ are:

\begin{itemize}
\item $\Delta < 0$:
$h(s) = \wp(s;g_2,g_3)$, $I=(0,2\omega_1)$.

\item $\Delta < 0$:
$h(s) = \wp_3(s;g_2,g_3) = \wp (s + \omega_3;g_2,g_3)$, $I=\R$.

\item $\Delta > 0$:
$h(s) = \wp(s;g_2,g_3)$, $I=(0,2\omega_1)$.

\item $\Delta =0$, $g_3 = -8a^3 > 0$:
\[
 h(s) = - 3a\tan^2{\left({\sqrt{-3a}}{s}\right)} -2a,
\quad I=(-\frac{\pi}{\sqrt{-12a}},\frac{\pi}{\sqrt{-12a}}).
     \]
\item $\Delta =0$, $g_3 = -8a^3 < 0$:
\[
 h(s) = 3a\tanh^2{\left({\sqrt{3a}}{s}\right)} -2a,
 \quad I=\R.
     \]

\item $g_2=g_3 =0$:
$h(s) = {s^{-2}}$, $I=(-\infty,0)$ or $I=(0,+\infty)$.

\end{itemize}

\section{Integration of the extremal trajectories}\label{s:integration}

For a potential function $h$ with invariants $g_2$, $g_3$, let
\[
 \mu_\pm (m,h) : = \frac{1}{2}\sqrt{P\left(\frac{m}{3} \pm 1;g_2,g_3\right)},
  \]
and, for each $s_0\in I$, define $\phi_\pm(m,h,s_0) : I \to \R$ by
\[
\phi_\pm(m,h,s_0) =
\begin{cases}
 \displaystyle\int_{s_0}^s{\frac{\mu_\pm(m,h)}{h(u) - \left(\frac{m}{3} \pm 1\right)}du},
  \quad \mu_\pm(m,h)\neq 0,\\
   \displaystyle\int_{s_0}^s{\frac{1}{h(u) - \left(\frac{m}{3} \pm 1\right)}du},
    \quad \mu_\pm(m,h)= 0.
     \end{cases}
       \]
Next, let $w_\pm(m,h)$ be the unique points\footnote{If $m=\pm 3$
and $g_2=g_3=0$, $w_\pm=\infty$.} in $\mathcal{D}(g_2,g_3)$ such
that
\[
  h(w_\pm) = \frac{m}{3} \pm 1 \quad\text{and} \quad   h'(w_\pm) = 2\mu_\pm(m,h).
   \]

Then, denoting by $\sigma_h$ and $\zeta_h$, respectively, the sigma
and zeta Weierstrassian functions corresponding to the potential
$h$,\footnote{$\sigma_h$ and $\zeta_h$ are the unique analytic odd
functions whose meromorphic extensions satisfy $\zeta'_h = -h$ and
$\sigma_h'/\sigma_h = \zeta_h$} we compute:

\vskip0.2cm
\noindent \textbf{Case I:} if $\mu_\pm(m,h) \neq 0$,
\[
 \phi_\pm(m,h,s_0) = \log{\frac{\sigma_h(s-w_\pm)}{\sigma_h(s+w_\pm)}} +
   2s\zeta_h(w_\pm) + c(s_0).
    \]

\vskip0.2cm
\noindent \textbf{Case II:} if $\mu_\pm(m,h) = 0$ and $g_2^2 +g_3^2 \neq 0$,
\[
 \phi_\pm(m,h,s_0) =\frac{-1}{3\left(\frac{m}{3} \pm 1\right)^2 -g_2/4}\left\{\zeta_h(s + w_\pm) +
  \left(\frac{m}{3} \pm 1\right)s\right\} + c(s_0).
   \]

\vskip0.2cm
\noindent \textbf{Case III:} if $\mu_\pm(m,h) = 0$ and $g_2=g_3= 0$,
\[
 \phi_\pm(m,h,s_0) =\frac{1}{3}s^3 + c(s_0).
   \]
Accordingly, define the maps
$R_\pm(m,h)$, $D_\pm(m,h,s_0) : I \to \GL(2,\C)$ as follows:

\vskip0.2cm
\noindent \textbf{Case I:} if $\mu_\pm(m,h) \neq 0$,
\[
 R_\pm(m.h) =\displaystyle \frac{1}{{2\mu_\pm}}
     \left(
     \begin{array}{cc}
     -\displaystyle\frac{{h'}-2\mu_\pm(m,h)}{2\sqrt{h - \left(\frac{m}{3} \pm 1\right)}}&
       -\sqrt{h - \left(\frac{m}{3} \pm 1\right)}\\
       \displaystyle\frac{{h'}+2\mu_\pm(m,h)}{2\sqrt{h - \left(\frac{m}{3} \pm 1\right)}}  &
         \sqrt{h - \left(\frac{m}{3} \pm 1\right)}\\
          \end{array}\right),
            \]
\[
D_\pm(m,h,s_0) =\left(
  \begin{array}{cc}
   \displaystyle \exp{(-\phi_\pm(m,h,s_0))}&
     \displaystyle 0\\
     \displaystyle 0 &
      \displaystyle \exp{(\phi_\pm(m,h,s_0))}\\
        \end{array}\right).
          \]

\vskip0.2cm
\noindent \textbf{Case II:} if $\mu_\pm(m,h) = 0$,
\[
 R_\pm(m.h) =
     \left(
     \begin{array}{cc}
     \displaystyle 1& 0\\
       \displaystyle \frac{{h'}}{2\sqrt{h - \left(\frac{m}{3} \pm 1\right)}} &
         \sqrt{h - \left(\frac{m}{3} \pm 1\right)}\\
          \end{array}\right),
            \]
\[
D_\pm(m,h,s_0) =\left(
  \begin{array}{cc}
   \displaystyle\frac{1}{\sqrt{h - (\frac{m}{3} \pm 1)}}&
     \displaystyle \phi_\pm(m,h,s_0)\\
     \displaystyle 0 &
      \displaystyle 1\\
        \end{array}\right).
          \]

Finally, define the maps $\Gamma_\pm(m,h,s_0)$, $\gamma(m,h,s_0) :
I\to \GL(2,\C)$ by
\[
 \Gamma_\pm(m,h,s_0) :=  R_\pm(m,h)(s_0)^{-1} D_\pm(m,h,s_0)(s_0)^{-1}
   D_\pm(m,h,s_0) R_\pm(m,h),
   \]
\[
 \gamma(m,h,s_0) :=  \Gamma_+(m,h,s_0)\, \Gamma_-(m,h,s_0)^{-1}.
   \]

\vskip0.1cm We are now in a position to state the following.

\begin{thm}\label{A}
The curve $\gamma(m,h,s_0)$ takes values in $\H^3_1 \cong \SL(2,\R)$
and defines an extremal trajectory with multiplier $m$ and reduced curvature
$h$. In particular, any extremal trajectory is equivalent to a curve of this
type.
\end{thm}

\begin{proof}
A direct, lengthy computation shows that
\begin{equation}\label{gamma-acca}
  \Gamma_\pm(m,h,s_0)^{-1}\, \Gamma_\pm(m,h,s_0)' = H(m,h)_\pm,
   \end{equation}
where
\begin{equation}\label{mh-hamiltoniani}
  H_+(m,h) = \left(
   \begin{array}{cc}
     0&1\\
      2h + \frac{m}{3} + 1&0\\
       \end{array}\right),\quad
  H_-(m,h) = \left(
   \begin{array}{cc}
     0&1\\
      2h + \frac{m}{3} - 1&0\\
      \end{array}\right).
       \end{equation}
Equations \eqref{gamma-acca} and \eqref{mh-hamiltoniani} imply that
$\Gamma_\pm(m,h,s_0)$ take values in a left coset of $\SL(2,\R)$.
Our normalization implies that $\Gamma_\pm(m,h,s_0)(s_0) = I_{2}$,
and hence that $\Gamma_\pm(m,h,s_0)$ and $\gamma(m,h,s_0)$ take
values in $\SL(2,\R)$. On the other hand, \eqref{gamma-acca} and
\eqref{mh-hamiltoniani} imply that $\Gamma(m,h,s_0)
=\left(\Gamma_+(m,h,s_0),\Gamma_-(m,h,s_0)\right)$ is a spinor frame
field with multiplier $m$ and curvature function $k = 2h
+\frac{m}{3}$. Therefore, $\gamma(m,h,s_0)$ is a trajectory with
curvature function $k$. This proves the required result.
\end{proof}

\section{Periodic trajectories}\label{s:periodic-solutions}

A trajectory is said \textit{quasi-periodic} if its reduced
curvature is a periodic function and if $\mu_\pm(m,h)$ are purely
imaginary. Let $Q\left(t; \ell,e_1\right)$ be the cubic polynomial
\[
 Q\left(t; \ell,e_1\right) = 4t^3 - 12\frac{(\ell^4 - \ell^2 +1)e_1}{(2-\ell^2)^2}t
  +4 \frac{(2\ell^4 + \ell^2 -1)e^3_1}{(2-\ell^2)^2}.
   \]
The reduced curvature of a quasi-periodic trajectory can be written in the form
\begin{equation}\label{potenziale-quasi-per}
 h(s;\ell,e_1) = e_1 \left(\frac{3\ell^2}{2-\ell^2}\mathrm{sn}^2\left(
\sqrt{\frac{3e_1}{2-\ell^2}}s,\ell\right) - \frac{1 +\ell^2}{2-\ell^2} \right),
   \end{equation}
where $(m,\ell,e_1)$ belongs to
$$
 \mathcal{W} = \{(m,\ell,e_1)\in\R^3 \,|\,
  e_1>0,\, \ell\in (0,1),\, Q\left(\frac{m}{3}\pm 1; l,e_1\right)<0\}.
  $$
The period of $h(s;\ell,e_1)$ is given by
\[
 p(\ell,e_1) = 2\sqrt{\frac{2-\ell^2}{3e_1}}K(\ell),
  \]
where $K$ is the complete elliptic integral of the third kind. We put
\[
 \rho_\pm(m,\ell,e_1) = \frac{1}{2}\sqrt{-Q\left(\frac{m}{3} \pm 1;\ell,e_1\right)}
  \]
and let $\mathcal{P} = (\mathcal{P}_+,\mathcal{P}_-) : \R \times \mathcal{W} \to
\R^2$ be the analytic map
\[
 \mathcal{P}_\pm(s,m,\ell,e_1) = \displaystyle\int_0^1{\frac{\rho_\pm(m,\ell,e_1)}
  {h(u;m,\ell,e_1) -\left(\frac{m}{3} \pm 1\right)}}du.
   \]
Recall that the spinor frame fields of $h$ are of the form
\[
\Gamma_\pm(s) = C\cdot \left(
  \begin{array}{cc}
   \displaystyle \exp{(i\mathcal{P}_\pm(s,m,\ell,e_1))}&
     \displaystyle 0\\
     \displaystyle 0 &
      \displaystyle \exp{(-i\mathcal{P}_\pm(s,m,\ell,e_1))}\\
        \end{array}\right)
        \cdot M(s),
          \]
where $C \in \GL(2,\C)$ and $M : \R \to \GL(2,\C)$ is a periodic map
with period $p(\ell,e_1)$. If we set $\Pi(m,\ell,e_1) =
\mathcal{P}(p(\ell,e_1),m,\ell,e_1)$, for each $(m,\ell,e_1) \in
\mathcal{W}$, then a quasi-periodic trajectory with invariants
$(m,\ell,e_1)\in \mathcal{W}$ is periodic if and only if
$\Pi(m,\ell,e_1) \in \mathbb{Q}\times \mathbb{Q}$.

We can state the following.

\begin{thm}\label{B}
There exists a discrete subset $\mathcal{M}_0\subset \R$ such that,
for every $m\in \R \setminus \mathcal{M}_0$, there exist countably
many periodic trajectories with multiplier $m$.
\end{thm}

\begin{figure}[ht]
\includegraphics{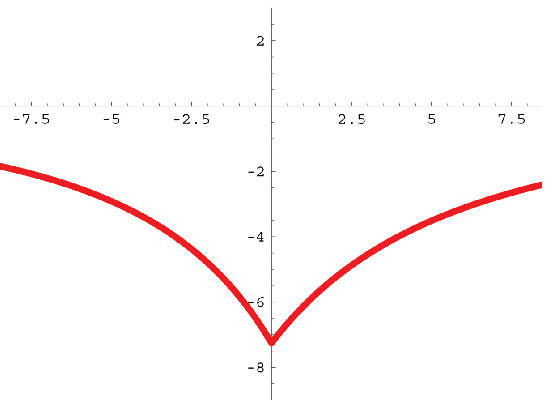}
\caption{The graph of the function $f$.}
\end{figure}

\begin{proof}
Consider the analytic map
\[
 \Psi : \mathcal{W} \ni (m,\ell,e_1) \mapsto
  \frac{\partial(\Pi_+,\Pi_-)}{\partial(\ell,e_1)}\vert_{(m,\ell,e_1)} \in \R.
   \]
If $\Psi(m,\tilde{\ell},\tilde{e_1}) \neq 0$, then $\Pi_m :
(\ell,e_1) \to \Pi(m,\ell,e_1) \in \R^2$ is a local diffeomorphism
near $(\tilde{\ell},\tilde{e_1})$. In this case there exist
countably many closed trajectories with multiplier $m$. The mapping
$\Psi$ can be computed explicitly, or numerically. To avoid dealing
with quite long formulae, we would rather adopt the numerical
viewpoint. Once we know $\Psi$, we define $f :  \R \to \R$, $m
\mapsto 400\Psi(m,\frac{1}{4}, |m| +10)$. This is an analytic
function for $m<0$ and for $m>0$. Looking at the graph of $f$ (cf.
Figure 1), we see that $f \neq \text{const.}$ This implies that the
set $\mathcal{M}_0 =\{m\in \R \, | \, f(m) =0\}$ is a discrete set
and that, for every $m \in \R \setminus \mathcal{M}_0$, there exist
countably many closed trajectories with multiplier $m$.
\end{proof}

%%%%%%%%%%%%%%%%
%\begin{figure}[ht]
%\includegraphics[scale=1.5]{Programma.eps}
%\caption{Program to compute the functions $\Psi$ and $f$ in Theorem
%2 (\textsc{Mathematica} 5.1)}
%\end{figure}
%%%%%%%%%%%%%%%

\section{Appendix: Derivation of the Euler--Lagrange equation}\label{s:appendix}

In this section, we outline the derivation of the Euler-Lagrange
equation associated with the $G$-invariant functional
\eqref{functional}. We follow a general construction for invariant
variational problems with one independent variable (due to
Griffiths) and write the Euler--Lagrange equation as a Pfaffian
differential system (PDF) $(\mathcal{J},\omega)$ on an associated
manifold $Y$. We adhere to the terminology and notations used in
\cite{Gr}.

\subsection{The variational problem}

The starting point of the construction is the replacement of the
original variational problem on null curves in $\SL (2,\R)$ by a
$G$-invariant variational problem for integral curves of a Pfaffian
differential system (PDF) with an independence condition on $M : =
G\times \R$.

Let $(\alpha, \beta)$ be the Maurer--Cartan form of $G = \SL (2,\R)
\times \SL (2,\R)$, where
\begin{equation}
 \alpha =
  \left(
   \begin{array}{cc}
     \alpha^1_1&\alpha^1_2\\
      \alpha^2_1&-\alpha^1_1\\
       \end{array}\right),
\quad \beta = \left(
   \begin{array}{cc}
     \beta^1_1&\beta^1_2\\
      \beta^2_1&-\beta^1_1\\
       \end{array}\right).
       \end{equation}
The Maurer-Cartan equations of $G$, or the structure equations, are
given by
\[
\begin{cases}
 d\alpha^1_1 = -\alpha^1_2 \wedge \alpha^2_1\\
 d\alpha^2_1 = 2\alpha^1_1 \wedge\alpha^2_1 \\
 d\alpha^1_2 = -2\alpha^1_1 \wedge\alpha^1_2
 \end{cases}\quad
\begin{cases}
 d\beta^1_1 = -\beta^1_2 \wedge\alpha^2_1 \\
 d\beta^2_1 = 2\beta^1_1 \wedge\alpha^2_1 \\
 d\beta^1_2 = -2\beta^1_1 \wedge\alpha^1_2.
 \end{cases}
   \]

On $M := G \times \R$, consider the PDS $(\mathcal{I}, \omega)$
defined by the differential ideal $\mathcal{I}$ generated  by the
linearly independent 1-forms
\[
 \begin{cases}
  \eta^1 = \frac{1}{2}(\alpha^1_2 -\beta^1_2), \quad
   \eta^2 = \alpha^1_1, \quad
    \eta^3 = \beta^1_1,\\
    \eta^4 = \frac{1}{2}(\alpha^1_2 + \beta^1_2) - \omega,\quad
     \eta^5 = \frac{1}{2}(\alpha^2_1 +\beta^2_1) -k \omega,
     \end{cases}
      \]
where
\[
 \omega := \frac{1}{2}(\alpha^2_1 -\beta^2_1)
  \]
gives the independence condition $\omega \neq 0$. Now, if $\gamma :
I \to \SL(2,\R)$ is a null curve without flex points, then the curve
$(\Gamma_\gamma;k_\gamma) : I \to M$, whose components are,
respectively, the spinor frame field $\Gamma = (\Gamma_+, \Gamma_-)$
along $\gamma$ and the curvature of $\gamma$, is an integral curve
of $(\mathcal{I}, \omega)$ (cf. Section \ref{s:pre}). Conversely,
any integral curve $(\Gamma,k) : I \to M$ of $(\mathcal{I}, \omega)$
defines a null curve with no flex points $\gamma : I \to \SL(2,\R)$,
where $\Gamma$ is the spinor field along $\gamma$, and $k$ is the
curvature of $\gamma$. So, null curves without flex points in
$\SL(2,\R)$ are identified with the integral curves of
$(\mathcal{I}, \omega)$.

From the above discussion, it follows that a null curve $\gamma$
without flex points is an extremal trajectory if and only if the
pair $(\Gamma_\gamma, k_\gamma)$ of its spinor frame field
$\Gamma_\gamma$ and curvature $k_\gamma$ is a critical point of the
functional defined on the space $\mathcal{V}(\mathcal{I},\omega)$ of
integral curves of $(\mathcal{I},\omega)$ by
\begin{equation}\label{var-action}
 \widehat{\mathcal{L}}: (\Gamma,k) \in \mathcal{V}(\mathcal{I},\omega)
  \mapsto \int_{I_{(\Gamma,k)}}{(\Gamma,k)^\ast ((m+k)\omega)}
   \end{equation}
when one considers compactly supported variations. Here,
$I_{(\Gamma,k)}$ is the domain of definition of $(\Gamma,k)$.

\subsection{The Euler--Lagrange system}

Following Griffiths \cite{Gr}, the next step is to associate to the
variational problem \eqref{var-action} a PDS $(\mathcal{J},\omega)$
on a new manifold $Y$, whose integral curves are stationary for the
associated functional.

For this, let $Z\subset T^\ast M$ be the affine subbundle defined by
\[
  Z = (m+k)\omega + I \subset T^\ast M,
   \]
where $I$ is the subbundle of $T^\ast M$ associated to the
differential ideal $\mathcal{I}$. The 1-forms $(\eta^1$, $\dots$,
$\eta^5$, $\omega)$ induce a global affine trivialization of $Z$,
which may be identified with $M\times \R^5$ by setting
\[
 M\times \R^5 \ni ((\Gamma,k); x_1,\dots,x_5)  \mapsto
 \omega_{|(\Gamma,k)} +{x_j\eta^j}_{|(\Gamma,k)} \in Z
 \]
(we use summation convention). Accordingly, the Liouville
(canonical) 1-form of $T^\ast M$ restricted to $Z$ is given by
\[
 \mu = (m + k)\omega + x_j\eta^j.
  \]
Exterior differentiation and use of the structure equations give
\[
\begin{split}
d\mu &\equiv dk \wedge \omega + (m +k)\left[(1+k)\eta^2 +
(1-k)\eta^3\right]
\wedge\omega + dx_j\wedge\eta^j\\
 &\quad + x_1(\eta^3 - \eta^2)\wedge\omega + x_2
\left[\eta^5 - (1+k)(\eta^1+\eta^4)\right] \wedge\omega\\
 &\quad + x_3\left[\eta^5 + (k-1)(\eta^1-\eta^4)\right]\wedge\omega
 - x_4 \left[(2+k)\eta^2 + (2-k)\eta^3\right] \wedge\omega \\
 &\quad +x_5\left[(1-k^2)\eta^2 - (1-k^2)\eta^3)\right]\wedge \omega
 - x_5 dk \wedge \omega
 \quad \mod \{\eta^i\wedge\eta^j\}.
 \end{split}
  \]

\vskip0.2cm

Then, we compute the Cartan system $\mathcal{C}(d\mu) \subset
T^{\ast}Z$ determined by the 2-form $d\mu$, i.e., the PDF generated
by the 1-forms $\left\{ i_\xi d\mu \, |\, \xi \in \mathfrak{X}(Z)
\right\} \subset \Omega^1(Z)$. Contracting $d\mu$ with the vector
fields of the tangent frame $
 \left(\frac{\partial}{\partial{\omega}},\frac{\partial}{\partial{k}},
  \frac{\partial}{\partial{\eta^i}},\frac{\partial}{\partial{x_j}}\right)
      $
on $Z$, dual to the coframe $\left(\omega, dk, \eta^i,dx_j\right)$,
$i,j =1,\dots,5$, we establish the following.

\begin{lemma}
The Cartan system $\mathcal{C}(d\mu)$, with independence condition
$\omega$, is generated by the 1-forms $\eta^1,\dots \eta^5$ and
\[
\begin{cases}
 \pi_1 = (1-x_5)dk, \quad \pi_2 = (1-x_5)\omega,  \quad
  \beta_1  = dx_1 +\left[x_2(1+k) +x_3(1-k)\right]\omega,  \\
   \beta_2 = dx_2 + \left[x_1 - (m+k)(1+k)+x_4(k+2) +x_5(k^2-1)\right]\omega, \\
    \beta_3 = dx_3 + \left[-x_1 - (m+k)(1-k)+x_4(2-k) +x_5(1-k^2)\right]\omega,\\
   \beta_4  = dx_4+ \left[x_2(1+k) +x_3(k-1)\right]\omega, \quad
    \beta_5 =  dx_5 -(x_2 +x_3)\omega.
     \end{cases}
      \]
\end{lemma}

The \textit{Euler--Lagrange system} associated to the variational
problem is the PDS $(\mathcal{J},\omega)$ on a submanifold $Y
\subset Z$ obtained by computing the involutive prolongation of
$(\mathcal{C}(d\mu),\omega)$. The submanifold $Y$ is called the
\textit{momentum space}. A direct calculation gives the following.

\begin{lemma}
The momentum space $Y\subset Z$ is defined by the equations
\[
 x_5 = 1, \quad x_2 +x_3 =0, \quad x_4 = \frac{m+k}{2}.
  \]
The Euler--Lagrange system $(\mathcal{J},\omega)$ is the PDF on $Y$
with independence condition $\omega$ generated by the 1-forms
${\eta^1}_{|Y},\dots,{\eta^5}_{|Y}$ and
\[
\begin{cases}
% {\eta^1}_{|Y},\dots,{\eta^5}_{|Y},\quad
  \sigma_1 =  dx_1 + 2kx_2\omega,\\
  \sigma_2 =  dx_2 + \displaystyle\left(\frac{k^2}{2} -\frac{mk}{2} -1 + x_1 \right)\omega,\\
   \sigma_3=  dk + 4 x_2 \omega.\\
    \end{cases}
       \]
\end{lemma}

\begin{remark}
The importance of this construction is that the projection $\pi_Y :
Y \to M$ maps integral curves of the Euler-Lagrange system to
extremals of the variational problem associated to
$(M,\mathcal{I})$. In our case, the converse is also true (see also
below), so that all extremals arise as projections of integral
curves of the Euler--Lagrange system. The theoretical reason for
this is that all derived systems of $\mathcal{J}$ have constant rank
(cf. \cite{Br}, \cite{Hsu}).

A direct calculation shows that $\mu_{|Y} \wedge (d\mu_{|Y})^4 \neq
0$ on $Y$, i.e., the variational problem is
nondegenerate.\footnote{A variational problem is said to be
\textit{nondegenerate} in case
\[
 \dim Y =2m+1 \quad \text{and}\quad
   \mu_{|Y} \wedge (d\mu_{|Y})^m \neq 0.
    \]
} This implies that $\mu_{|Y}$ is a contact form and that there
exists a unique vector field $\zeta \in \mathfrak{X}(Y)$, the
\textit{characteristic vector field} of the contact structure, such
that $\mu_{|Y}(\zeta) =1$ and $i_\zeta\,d\mu_{|Y} =0$. In
particular, the integral curves of the Euler-Lagrange system
coincide with the characteristic curves of $\zeta$.

\end{remark}

\subsection{The Euler--Lagrange equation}\label{ss:e-l-eq}

Let $\mathcal{V}(\mathcal{J},\omega)$ be the set of integral curves
of the Euler-Lagrange system. If $y=((\Gamma,k);x_1,x_2) : I \to Y$
is in $\mathcal{V}(\mathcal{J},\omega)$, the equations
\[
 \eta^1=\eta^2 =\cdots =\eta^5 =0
  \]
and the independence condition $\omega \neq 0$ tell us that
$\Gamma=(\Gamma_+,\Gamma_-)$ defines a spinor frame along the null
curve $\gamma = \Gamma_+ \Gamma_-^{-1}$ and that $k$ is the
curvature of $\gamma$.

Next, for the smooth function $k : I \to \R$, define $k'$, $k''$ and
$k'''$ by
\[
 dk = k'\omega,\quad
  dk' = k''\omega, \quad dk'' = k''' \omega.
    \]
Equation $\sigma_3=0$ implies
\[
 x_2 = -\frac{k'}{4}.
  \]
Further, equation $\sigma_2=0$ gives
\[
 x_1 = \frac{k''}{4}-\frac{k^2}{2} + \frac{mk}{2} +1.
  \]
Finally, equation $\sigma_1=0$ yields
\begin{equation}\label{euler-lagrange-eq}
  k''' -6 k k' + 2  m k' = 0.
   \end{equation}
This coincides with the Euler--Lagrange equation of the extremals of
\eqref{functional}, which has been computed in \cite{FGL-PhysLettB}.
Thus, an integral curve of the Euler--Lagrange system projects to an
extremal trajectory in $\SL(2,\R)$.

Conversely, if $\gamma : I \to \SL(2,\R)$ is a null curve without
flex points, $\Gamma_\gamma$ its spinor frame, and $k_\gamma$ its
curvature, let $y_\gamma : I \to Y$ be the lift of $\gamma$ to $Y$
given by
\[
  y_\gamma(t) =
\left((\Gamma_\gamma, k_\gamma); \frac{k''}{4}-\frac{k^2}{2} +
\frac{mk}{2} +1,
   -\frac{k'}{4}\right).
    \]
Then, $y_\gamma$ is an integral curve of the Euler--Lagrange system
if and only if $k_\gamma$ satisfies equation
\eqref{euler-lagrange-eq} if and only if $\gamma$ is an extremal
trajectory. Thus, the integral curves of the Euler--Lagrange system
arise as lifts of trajectories in $\SL(2,\R)$.

\begin{remark}
Griffiths' approach to calculus of variations, besides for providing
the Euler--Lagrange equations, is important for giving an effective
procedure to construct the momentum mapping induced by the
Hamiltonian action of $G$ on $Y$ and to prove that it is constant on
the integral curves of $(\mathcal{J},\omega)$, which in turn leads
to the integration by quadratures of the extrema (cf. \cite{GM},
\cite{MN-SICON}).
\end{remark}

\bibliographystyle{amsalpha}

\begin{thebibliography}{AA}

\bibitem{Br}
R. L. Bryant, {On notions of equivalence of variational problems
with one independent variable}, \textit{Contemp. Math.} \textbf{68}
(1987), 65--76.

\bibitem{Br1987}
R.L. Bryant, Surfaces of mean curvature one in hyperbolic space,
Th\'eorie des vari\'et\'es minimales et applications (Palaiseau,
1983--1984), Ast\'erisque, No. \textbf{154-155} (1987), 12,
321--347, (1988).

\bibitem{BFJL-CQG}
M. Barros, A. Ferr\'andez, M.A. Javaloyes, P. Lucas, Relativistic
particles with rigidity and torsion in $D=3$ spacetimes,
\textit{Classical Quantum Gravity} \textbf{22} (2005), no. 3,
489--513.

\bibitem{FGL-PhysLettB}
A. Ferr\'andez, A. Gim\'enez, P. Lucas, Geometrical particle models
on 3D null curves, \textit{Phys. Lett. B} \textbf{543} (2002),
311--317; hep-th/0205284.

\bibitem{GMM-MA}
J.A. G\'alvez, A. Mart\'{\i}nez, F. Mil\'an, Flat surfaces in the
hyperbolic 3-space, \textit{Math. Ann.} \textbf{316} (2000),
419--435.

\bibitem{GM}
J.D.E. Grant, E. Musso, Coisotropic variational problems, \textit{J.
Geom. Phys.} \textbf{50} (2004), 303--338; math.DG/0307216.

\bibitem{Gr}
 P. A. Griffiths,
\textit{Exterior differential systems and the calculus of
variations}, Progr. Math., 25, Birkh\"auser, Boston, 1982.

\bibitem{Hsu}
L. Hsu, {Calculus of variations via the Griffiths formalism},
\textit{J. Differential Geom.} \textbf{36} (1992), 551--589.

\bibitem{Inoguchi}
J. Inoguchi, S. Lee,
Null curves in Minkowski 3-space, preprint, 2006.

\bibitem{MN-SICON}
E. Musso, L. Nicolodi, Reduction for constrained variational
problems on 3D null curves, submitted to \textit{SIAM J. Control
Optim.}.

\bibitem{NMMK-NucPhysB}
A. Nersessian, R. Manvelyan, H.J.W. M\"uller-Kirsten, Particle with
torsion on 3d null-curves, \textit{Nuclear Phys. B} \textbf{88}
(2000), 381--384; hep-th/9912061.

\bibitem{Pi}
E. Picard, {Sur les \'equations diff\'erentielles lin\'eaires \`a
coefficients dou\-ble\-ment p\'eriodiques}, \textit{J. Reine Angew.
Math.}
%(Crelle's Journal)
\textbf{90} (1881), 281--302.

\bibitem{UY-Crelle}
M. Umehara, K. Yamada, A parametrization of the Weierstrass formulae
and perturbation of  complete minimal surfaces in $\R^3$ into the
hyperbolic  3-space, \textit{J. Reine Angew. Math.} \textbf{432}
(1992), 93--116.

\end{thebibliography}

\end{document}